\documentclass{amsart}

\usepackage[mathcal]{euscript}
\usepackage{mathrsfs}
\usepackage{amsfonts,amssymb,amsbsy,amsmath}
\usepackage{graphicx,color}
\usepackage{epsfig}
\usepackage[all,cmtip]{xy}

\newcommand{\N}{\mathbb{N}}
\newcommand{\D}{\mathbb{D}}
\newcommand{\F}{{\mathbb{F}}}

\newcommand{\CP}{{\mathbb{C}}P}
\newcommand{\R}{\mathbb{R}}
\newcommand{\Q}{\mathbb{Q}}
\newcommand{\C}{\mathbb{C}}
\newcommand{\Z}{\mathbb{Z}}
\newcommand{\M}{\mathit{Map^+}}
\newcommand{\Mm}{\mathit{Map}}

\newcommand{\ba}{\begin{array}}
\newcommand{\ea}{\end{array}}

\newcommand{\Sg}{\Sigma_g}

\theoremstyle{plain}
\newtheorem{t.}{Theorem}
\newtheorem{l.}[t.]{Lemma}
\newtheorem{p.}[t.]{Proposition}
\newtheorem{c.}[t.]{Corollary}

\theoremstyle{definition}
\newtheorem*{d.}{Definition}

\newtheorem*{r.}{Remark}

\begin{document}

\title{Real elements in the mapping class group of  $T^2$}

\author{Nerm\.{\i}n Salepc\.{\i}}
\address{Universit\'e de Pau et des Pays de l'Adour
Avenue de l'Universit\'e - BP 1155
Batiment IPRA 
64013 PAU CEDEX}
\email{nermin.salepci@univ-pau.fr}

\begin{abstract}
We present a complete classification of elements in the mapping class group of the torus which have a representative that can be written as a product of two orientation reversing involutions. Our interest in such decompositions is motivated by  features of the monodromy maps of real fibrations. We employ the property that the mapping class group of the torus is identifiable with $SL(2,\Z)$ as well as that the quotient group $PSL(2,\Z)$ is the symmetry group of the {\em Farey tessellation} of the Poincar\'e disk. 
\end{abstract}

\maketitle
\section{Introduction}
In this work, we are interested in finding the isotopy classes of orientation preserving diffeomorphisms of  the torus which have a representative that decomposes into a product of two orientation reversing involutions. The question arose when we were studying properties of the monodromy maps of real Lefschetz fibrations. 

Let us recall that a Lefschetz fibration is  a projection of an oriented smooth 4-manifold onto an oriented surface such that
apart from finitely many fibers, which have a single node, the fibers are actually smooth oriented surfaces. 
Intuitively, a real structure can be regarded as a topological generalization of the complex conjugation. We define a {\em real structure} on an oriented 4-manifold as an orientation preserving involution whose fixed point set, if it is not empty, has dimension 2.  Likewise, a {\em real structure} on an oriented surface is an orientation reversing involution. A smooth manifold together with a real structure is called a {\em real manifold} and the fixed point set of a real structure is called the {\em real part}.   
A {\em real Lefschetz fibration}  is, thus, defined as  a Lefschetz fibration of a real 4-manifold over a real surface where the fiber structure is compatible with the real structures.

The fundamental fact about real Lefschetz fibrations is that their monodromy maps along the loops on which the real structure induces an orientation reversing action decompose into a product of two orientation reversing involutions. These involutions are the real structures of the two real fibers over the two real points of the loop. The diffeomorphisms (as well as their classes) with such a feature  are called {\em real}.  Let us remark that if, in particular, the base space is a surface with a single boundary component, then  the monodromy map along the boundary component satisfies this property.

We can rephrase the decomposition property as follows: the monodromy maps as above of real Lefschetz fibrations are conjugated to their inverses by a real structure.  At this point, a weaker property, that of being conjugated to its inverse by an orientation reversing diffeomorphism  appears naturally. Diffeomorphisms (as well as their classes) with the latter property are called {\em weakly real}.
We study weakly real and real diffeomorphisms simultaneously.

As is well known, the restriction of a Lefschetz fibration to  loops  in the complement of the critical set is a usual fibration over a circle with fiber compact connected oriented smooth $2$-manifold $F$.  We call such fibrations  {\em $F$-fibrations} and extend the above definitions to $F$-fibrations to proceed with them.

The aim of this article is to answer the following questions in the case of $F=T^2$ (elliptic):

\begin{enumerate}

\item  Which classes are real/ weakly real?  

\item  Obviously real classes are weakly real. What about the converse:  are weakly real classes real? 

\end{enumerate}

Employing the identification of the mapping class group of the torus with  $SL(2,\Z)$ we state the main results as follows:

\begin{t.} 
All elliptic and parabolic matrices in  $SL(2,\Z)$ are real.

A hyperbolic matrix $A \in SL(2,\Z)$ is real if and only if its cutting period-cycle  $[a_1 a_2 \ldots a_{2n}]_A$ is odd-bipalindromic.

Moreover, a matrix in $SL(2,\Z)$ is real if and only if it is weakly real.
\end{t.}

(The three statements of the above theorem are indeed presented separately as Theorem~\ref{elpar}, Theorem~\ref{hiper} and Theorem~\ref{wrr}, respectively.)

Elliptic, parabolic and hyperbolic matrices in $SL(2,\Z)$ differ by the nature of their fixed points.  In Section~\ref{eliparabol}, 
we discuss the conjugacy classes of elliptic and parabolic matrices.  Conjugacy classes of hyperbolic matrices are explained in 
Section~\ref{hiperbol}, where  we also  give the definition of the cutting period-cycle. The last two sections  are devoted to the real factorizations of the matrices. We define odd-bipalindromic cutting period-cycle in the last section.

Finally, we would like to note that this work is self-contained and the tools used are explained in detail.

\textbf{Acknowledgements.}
This work is a part of my PhD thesis where we introduced real Lefschetz fibrations and gave a classification of real elliptic Lefschetz fibrations over $S^2$ with only real critical values.  I would like to express my gratitude to my supervisors  Sergey Finashin and Viatcheslav Kharlamov  for  their guidance and endless support throughout  my research.  

I would also like to thank  Caroline Series for her interest to my questions and for sending me her articles and  Allen Hatcher for pointing out the reference I needed as well as Slava Matveyev for fruitful discussions.

\section{Real  $F$-fibrations and their monodromies}
Let  $Y$ be a compact connected oriented smooth 3-manifold,
and  $\pi:Y\to S^1$ be a fibration whose fiber is an
oriented smooth $2$-manifold $F$. We call
$\pi$ an {\em $F$-fibration}. In particular, when $F=T^2$, we call $\pi$ an \emph{elliptic $F$-fibration}.

\begin{d.} \label{realS1}
An $F$-fibration $\pi:Y\to S^1$ is called \emph{weakly real} if
there is an orientation preserving diffeomorphism $H:Y\to Y$ which
sends fibers into fibers reversing their orientations. If
$H^2=id$, then $H$ will be called a \emph{real structure} on the
$F$-fibration $Y\to S^1$. An $F$-fibration equipped with a real
structure will be called \emph{real}. 

Note that $H$ induces an orientation reversing diffeomorphism
$h:S^1\to S^1$ such that the following diagram commutes:

$$\xymatrix{Y \ar[r]^{H} \ar[d]_{\pi} & Y \ar[d]^{\pi} \\
          S^1 \ar[r]^{h}  & S^1.}$$
\end{d.}

Here, it is not difficult to see that the set of orientation reversing
involutions form a single conjugacy class in the diffeomorphism
group of $S^1$ (the crucial observation is that any such
involution has precisely two fixed points).
Therefore, any real $F$-fibration is equivariantly isomorphic
to an $F$-fibration whose involution $h$ is standard. Let it be the complex conjugation $c_{S^1}: S^1\to S^1$,
$z\mapsto\bar z$, $z\in S^1\subset\C$.

In the case of a weakly real $F$-fibration, $h$ may be not
an involution; however, it also has precisely two fixed points and
can be changed into an involution by an isotopy. This isotopy can be lifted to an isotopy of
$H$. Thus, by modification of $H$ we can always make $h$ an
involution.
 So it is not restrictive for us to suppose always that  $h=c_{S^1}$ for both
real and weakly real $F$-fibrations.

The restrictions of $H$ to the invariant fibers
$F_\pm=\pi^{-1}(\pm1)$ will be denoted $h_\pm:F_\pm\to F_\pm$.
 In the case of real $F$-fibrations, we will prefer using  the notation
$c_Y$ for the involution $H$, and $c_\pm$ for the involutions
$h_\pm$.

In what follows, we choose the point $b$ in the upper semi-circle,
$S^+$.  The restriction $Y_+=\pi^{-1}(S^+)\to S^+$ of $\pi$
admits a trivialization $\phi_+:Y_+\to F\times S^+$ which is
identical on the fiber $F=F_b$. In the case of real fibrations, this allows us to consider the
pull-back of $c_\pm$ via $\phi$, namely, the two involutions
$x\mapsto \phi_+(c_\pm(\phi_+^{-1}(x\times\pm1)))$ on the same
fiber $F$. We will stick to the notation $c_\pm$ for these
involutions.

It is well known that any
$F$-fibration $\pi:Y\to S^1$ is isomorphic to the projection
$M_f\to S^1$ of a mapping torus $M_f=F\times
I\diagup_{{(f(x),0)\sim (x,1)}}$ of some orientation preserving diffeomorphism $f:F\to
F$. More precisely, if we fix a particular fiber $F=F_b=\pi^{-1}(b)$, $b\in S^1$,
then an isomorphism $\phi:M_f\to Y$ can be chosen so that
$F\times0$ and $F\times1$ are identified with the fiber $F_b$, so
that $x\times0\mapsto x$ and $x\times1\mapsto f(x)$.

An $F$-fibration $\pi$ determines an orientation preserving diffeomorphism $f$
up to isotopy and thus, provides a well-defined element in the
mapping class group $[f]\in\M(F)$ called the {\em monodromy of
$\pi$} (relative to the fiber $F=F_b$).
A map $f$ representing the
class $[f]$ will be also often called {\em monodromy}, or more
precisely, a {\em monodromy map}.

In some cases, we fix a marking $\rho:\Sg\to F_b$ (an identification of $F_{b}$ with an abstract genus-$g$ surface $\Sigma_{g}$). Then the
diffeomorphism $\rho^{-1}\circ f\circ\rho:\Sg\to\Sg$ (the
pull-back of $f$) as well as its isotopy class $[\rho^{-1}\circ
f\circ\rho]\in\M(\Sg)$ will be called the {\em monodromy of $\pi$
relative to the marking $\rho$}.

\begin{t.}\label{fcc}
Let $\pi: Y \to S^1$ be a weakly real $F$-fibration with a distinguished fiber $F=F_b$, $b\in S^+$. Then the two
product diffeomorphisms of the fiber $F$, $(h_+)^{-1} \circ h_-$,
and $h_+\circ (h_-)^{-1}$ are isotopic and equal the monodromy
of $\pi$ relative to the fiber $F$.
In particular, if $\pi$ is a real $F$-fibration, then the monodromy can
be factorized as $c_+ \circ c_-$.
\end{t.}

\noindent {\it Proof:} 
Consider a trivialization $Y_-\to F\times S^-$ of the
restriction $Y_-=\pi^{-1}(S^-)\to S^-$ of $\pi$ over the lower
semi-circle, $S^-$, which is the composition of $\phi_+\circ H
:Y_-\to F\times S^+$, with the map $F\times S^+\to F\times S^-$,
$(x,z)\mapsto (x, c_{S^1}(z))$.

If $S^1$ is split into several arcs and a fibration over $S^1$ is
glued from trivial fibrations over these arc, then the monodromy
is clearly the product of the gluing maps of the fibers over the
common points of the arcs, ordered in the counter-clockwise
direction beginning from a marked point $b\in S^1$.
In our case, the arcs are $S^+$, $S^-$, their common points
follow in the order $-1$, $+1$, and the corresponding gluing maps,
are $h_-$ and $h_+^{-1}$. This gives monodromy $(h_+)^{-1}\circ
h_-$. If we consider another trivialization $Y_-\to F\times
S^-$ replacing in its definition $H$ by $H^{-1}$, then the
gluing maps will be $h_-^{-1}$ and $h_+$, and the monodromy is
factorized as $h_+\circ (h_-)^{-1} $.
\hfill  $\Box$\\ 

\begin{c.}\label{2cfc} 
(1) If $f:F \to F$ is a monodromy map of  a weakly real $F$-fibration, then the diffeomorphisms $h^{-1}\circ f \circ h$ as well as $h\circ f\circ h^{-1}$, where $h$ stands either for $h_+$, or for $h_-$, are all isotopic to the inverse $f^{-1}$.  

(2) If $f$ is a monodromy map  of a real $F$-fibration, then $f^{-1}=c_{+}\circ f\circ c_{+}=c_{-}\circ f\circ c_{-}$. \hfill  $\Box$\\ 
\end{c.} 

\begin{r.}  
It is obvious that  $f=c_{+} \circ c_{-}$  for some real structures $c_{-}, c_{+}$ if and only if  $f^{-1}=c_{\pm} \circ f \circ c_{\pm}$. It follows from the known cases of the Nielsen realization problem that  $[f]=[c_{+} \circ c_{-}] $  if and only if $[f^{-1}]=[c_{\pm} \circ f \circ c_{\pm}] $, (cf. \cite{He}, \cite{K}).  In another words,  if 
$[f^{-1}]=[c \circ f \circ c] $,  then there is a diffeomorphism $g$ isotopic to $f$ such that  $g^{-1}=c \circ g \circ c $.  Similarly, if  $[f^{-1}]=[h\circ f \circ h^{-1} ] $ for some orientation reversing diffeomorphism $h$, then there are diffeomorphisms $g, k$  such that $[g]=[f]$ and $[k]=[h]$ with  $g^{-1}=k\circ g \circ k^{-1}$.   (All  the class equalities above are considered  in the extended mapping class group $\Mm(T^2)$ of the torus.)   \end{r.}

\begin{d.}\label{zayifguclureel}
A diffeomorphism $f: F\to F$  as well as its isotopy class $[f]\in \M(F)$ will be called {\em real}  if there is a real structure $c:F \to F$ such that $f^{-1}=c\circ f\circ c$.

We call $f:F\to F$ as well as its isotopy class $[f]\in \M(F)$  {\em weakly real}  if $f^{-1}=h\circ f \circ h^{-1}$ for some orientation reversing diffeomorphism $h$ of $F$.  
\end{d.}

\begin{p.}\label{s1fcc}
An $F$-fibration is real (respectively weakly real) if and only if
its monodromy $f$ is real (respectively weakly real).
\end{p.}

\noindent {\it Proof:} 
We give the proof for real $F$-fibrations; the proof for weakly real ones is
analogous.

 Necessity of the condition follows from Corollary~\ref{2cfc}. 
As for the converse, consider an $F$-fibration $\pi:Y\to
S^1$ with the monodromy class $[f]\in \M(F)$,
and $f$ its representative such that   $f^{-1}=c\circ f\circ c$
for some real structure  $c$  on $F$. Presenting $Y$ as
$F\times [0,1]\diagup_{(f(x), 0)\sim(x,1)}$, we obtain a well-defined
involution $c_Y:Y\to Y$ induced from the involution
$(x,t)\mapsto(c(x),1-t)$ on $F\times [0,1]$.
 It preserves the fibration structure and acts as $c$ and $f\circ
 c$ on the real fibers $F\times\frac12$ and $F\times0=F\times1$
respectively.
 \hfill  $\Box$\\

\section{Homology monodromy factorization of elliptic $F$-fibrations}

It is well known that $\M(T^2)=SL(2,\Z)$, due to the fact that every
diffeomorphism $f:T^2 \rightarrow T^2$ is isotopic to a {\em linear
diffeomorphism}. The latter diffeomorphisms by definition are
induced on $T^2=\R^2/\Z^2$ by a linear map $\R^2\to\R^2$ defined by
a matrix $A\in SL(2,\Z)$. Note that we can naturally identify
$T^2=H_1(T^2,\R)/H_1(T^2, \Z)$ and interpret matrix $A$ as the induced
automorphism $f_*$ in $H_1(T^2, \Z)$. The latter automorphism is
 called the \emph{homology monodromy}. Since
isotopic diffeomorphisms have the same homology monodromy in
$H_1(T^2,\Z)$, we obtain well-defined homomorphisms $\M(T^2)\to
Aut^+(H_1(T^2,\Z))\to SL(2,\Z)$ which are in fact isomorphisms (here
$Aut^+$ stand for the orientation preserving automorphisms).
Let $a$ denote the simple closed curve on $T^2$ represented by the
equivalence class of the horizontal interval
$I\times0\subset\R^2$, and $b$ is similarly represented by the
vertical interval $0\times I$. We have $a\circ b=1$; hence, the
homology classes represented by these curves are integral
generators of $H_1(T^2,\Z)$. The mapping class group $\M(T^2)$ of $T^2$ is
generated by the Dehn twists $t_a$ and $t_b$, which can be
characterized by their homology monodromy homomorphism matrices
${t_a}_*=\tiny{\left(\begin{array}{cc}
           1 & 1\\
           0 & 1
         \end{array} \right)}$ and
${t_b}_*=\tiny{\left(\begin{array}{cc}
           1 & 0\\
           -1& 1
         \end{array} \right)}$.

In a like manner, it can be shown that  the extended mapping class group $\Mm(T^2)$  is isomorphic to the general linear group $GL(2,\Z)$. 
As a consequence, for elliptic $F$-fibrations, the question of
characterization of real monodromy classes $[f]\in\M(T^2)$ can be
interpreted as the question on the decomposability of their
homology monodromy $f_* \in SL(2,\Z)$ into a product of two {\em
linear real structures} in the group $GL(2,\Z)$. The latter structures by definition  are
linear orientation reversing maps of order 2 defined by integral
$(2\times 2)$-matrices.
 Such decomposability is equivalent to the
property that $f_*$ is conjugate to its inverse by a linear real
structure. Hence, a necessary condition for a matrix $A$ to be real
is that both $A$ and $A^{-1}$ lie in the same conjugacy classes
in the group $GL(2,\Z)$.

Recall that there are three types of real structures  on $T^2$
distinguished by the number of their real components: $0$, $1$, or
$2$. Note that the automorphisms of $H_1(T^2,\Z)$ induced by
the real structures with 0 or 2 real components  are diagonalizable over $\Z$, namely, their
matrices are conjugate to
 $\tiny{\left(\begin{array}{cc}
           1 & 0\\
           0 & -1
         \end{array} \right)}$
in $GL(2,\Z)$. So we cannot determine if the number of components
$0$ or $2$ knowing only the matrix representing the homology
action of the real structure. The homology action of a real structure with 1 real component is presented by a
matrix conjugate to
 $\tiny{\left(\begin{array}{cc}
           0 & 1\\
           1 & 0
         \end{array} \right)}$.

\section{The modular action on the hyperbolic half-plane}

Let $\C^2$ be considered as the vector space of $2\times1$ matrices over $\C$. Then a matrix $A=\tiny{\left(\begin{array}{cc}
           a & b\\
          c  & d
         \end{array} \right)}$ in $GL(2,\Z)$  acts on $\C^2$ from the left as matrix multiplication. $$\small{\left(\begin{array}{cc}
           a & b\\
          c  & d
         \end{array} \right)\left(\begin{array}{c}
           z_1\\
           z_2
         \end{array} \right)=\left(\begin{array}{c}
           az_1 + bz_2\\
          cz_1  + dz_2
         \end{array} \right)}$$

\begin{figure}[ht]
   \begin{center}
      \includegraphics[ width=6cm]{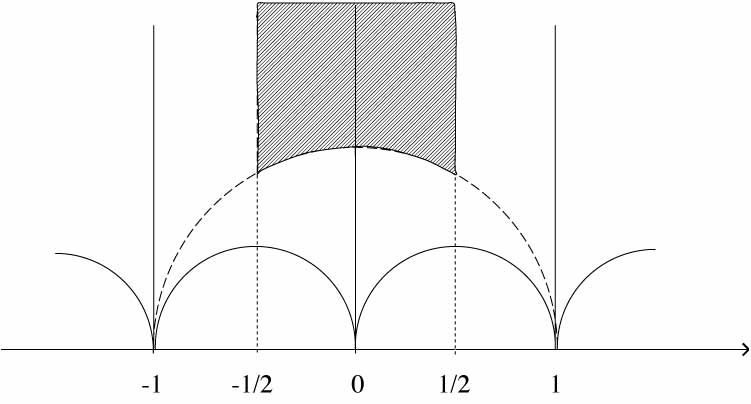}
           \caption{\small{The upper half plane model of the hyperbolic plane and the standard fundamental domain of $SL(2,\Z)$}.}
       \label{fundomain}
      \end{center}
 \end{figure}

This action can be extended to $\CP^1= \C^2 \setminus \{(0,0)\}\diagup_{\tiny{(z_1, z_2) \sim (\lambda z_1, \lambda z_2)}}$.
Let us identify $\CP^1\cong \{(z_1, z_2)\in \C^2, z_2\neq 0 \}\cup
\{\infty \} \cong  \C \cup \{\infty\}$  and rewrite the action of
$GL(2,\Z)$. We obtain a linear fractional transformation
$z\rightarrow  \frac{az+b}{cz+d}$ where $z=\frac{z_1}{z_2}$. In
particular, if $A \in SL(2 ,\Z)$, then the transformation preserves
the orientation of $\C$ and takes $\R \cup \{\infty\}$ to itself
preserving its orientation. Hence, it gives rise to a
diffeomorphism of the upper half plane $\mathbb{H}$ which can be
seen as a model for the hyperbolic plane where the geodesics are
the semi-circles centered at a real point or vertical half-lines
which can also be considered as arcs of infinite radius. By
identifying the upper half plane with the lower half plane by the  complex
conjugation, one extends the action of $SL(2,\Z)$ to an action of
$GL(2,\Z)$. The standard fundamental domain of the action is the
set $\{z| \left|Re(z)\right|\leq \frac{1}{2}, \left|z\right|\geq 1
\}$ which is shown in Fig.~\ref{fundomain}.

\section{The Farey Tessellation}

\begin{figure}[ht]
   \begin{center}
   \includegraphics[width=4.8cm, scale=0.3]{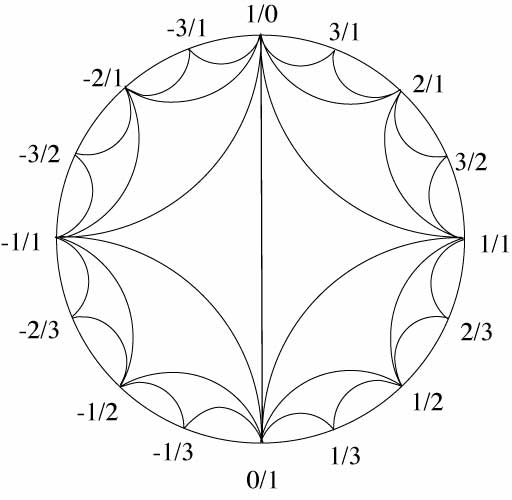}
   \caption{\small{Tessellation of $\D$.}}
       \label{tesel}
      \end{center}
 \end{figure}

Let us identify the upper half plane model with the Poincar\'e disk model $\D$. We will consider the disk $\D$ together with its boundary $\R \cup \infty$ and define a tessellation on $\D$ as follows:

Set $\infty$ as $\frac{1}{0}$ and consider the two fractions $\frac{0}{1}$ and $\frac{1}{0}$, spot them on $\D$ as the south and the north poles respectively and connect them with a line which will be the vertical diameter. Consider their mediant $\frac{0+1}{1+0}= \frac{1}{1}$ and connect each of them with a geodesic to the mediant. Apply the same to the fractions $\{\frac{0}{1}, \frac{1}{1}\}$ and $\{\frac{1}{1}, \frac{1}{0}\}$. Iterating this process one obtains a tessellation of the right semi-disk. By taking the symmetry one extends the tessellation to $\D$, see Fig.~\ref{tesel}.

In the literature this tessellation is called the \emph{Farey
tessellation}. Let us denote the disk
together with the Farey tessellation by $\D_F$. Note that the Farey
tessellation is a tessellation of $\D$ by ideal triangles ( i.e. triangles with vertices on the boundary $\D_F$).
In fact, the set of vertices of the triangles is exactly $\Q \cup \{\infty\}$.
Moreover, two fractions $\frac{m_1}{n_1},  \frac{m_2}{n_2}$ are
connected by a {\em line} if and only if $m_1n_2-m_2n_1=\pm1$. Hence, the action of
$GL(2,\Z)$ on $\D$ induces an action on $\D_F$ which is
transitive on the geodesics of $\D_F$. Only $\pm I$ acts as the identity; hence, the modular
group $PGL(2,\Z)=GL(2,\Z) / \pm I$ is the symmetry group of $\D_F$
where the subgroup $PSL(2,\Z) =SL(2,\Z) / \pm I$ gives the
orientation preserving symmetries. In what follows we denote by $\Gamma$ the triangle with vertices $\{0, 1,
\infty\}$. Note that $\Gamma$ splits in 3 copies of the fundamental region.

The fixed points of the modular action of a matrix $A\in
PSL(2,\Z)$, $A\neq I$, in $\D_F$ are solutions of
$z=\frac{az+b}{cz+d}$. This gives a quadratic equation
with the discriminant
$tr(A)^2-4$. If
the trace $\left| tr(A)\right|<2$, then the discriminant is
negative and the modular action is a rotation around an imaginary
point (an interior point of $\D_F$). Such matrices are called {\em
elliptic}. If $\left| tr(A)\right|=2$, then the discriminant
vanishes, and $A$ acts as a translation with one fixed rational
point, $\frac{d-a}2$ (on the boundary of $\D_F$). Such matrices
are called {\em parabolic}.
 The {\em hyperbolic matrices} have
$\left| tr(A)\right|>2$ and define a translation of $\D_F$ with
two fixed quadratically irrational real points on the boundary of
$\D_F$.

\section{Conjugacy classes of elliptic and parabolic matrices}\label{eliparabol}

\begin{figure}[ht]
   \begin{center}
      \includegraphics[width=4.8cm, scale=0.3]{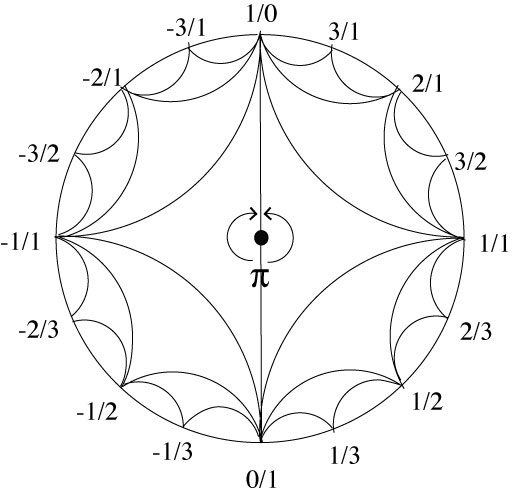}\hspace{.5cm}
   \includegraphics[width=4.8cm, scale=0.3]{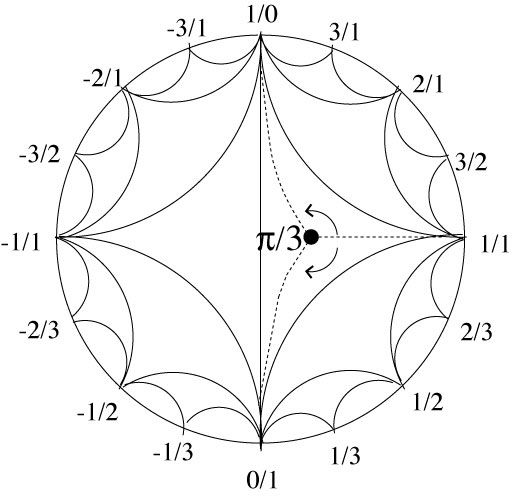}
       \caption{\small{Modular actions of elliptic matrices, $E_{\pm \pi}$,$\,E_{\pm \frac{2\pi}{3}}$.}}
       \label{elliptic}
      \end{center}
 \end{figure}

\textbf{Elliptic matrices:} As mentioned above an elliptic matrix  $A\in PSL(2,\Z)$ acts on $\D_F$ as a rotation around a point in the interior of $\D_F$. The center of the rotation belongs to one of the triangles of the tessellation. Without loss of generality, let us assume that the fixed point belongs to the triangle $\Gamma$.
If the fixed point belongs to an edge of $\Gamma$, then $A$ rotates $\Gamma$ by angles $\pm \pi$. Let us denote such rotations by $E_{\pm \pi}$.
The other possibility is the rotations $E_{\pm\frac{2\pi}3}$ by  angles $\pm\frac{2\pi}3$ around the center of $\Gamma$, see 
Fig.~\ref{elliptic}. It is not hard to see that the matrices  $E_{\pm \pi}$ (respectively $E_{\pm\frac{2\pi}3}$) are conjugate to each other via an orientation reversing matrix in $PGL(2,\Z)$.

Since $PGL(2,\Z)$ acts transitively on the triangles of the
tessellation, $E_\pi=\tiny{\left(\begin{array}{cc}
           0 & 1\\
          -1 & 0
         \end{array} \right)}$ and $E_{\frac{2\pi}{3}}=\tiny{\left(\begin{array}{cc}
           0 & 1\\
          -1 & 1
         \end{array} \right)}$ represent the conjugacy classes in $PGL(2,\Z) $ of elliptic matrices in $PSL(2,\Z)$.

Each matrix $A$ in $PSL(2,\Z)$ defines two matrices $\pm A$ in $SL(2,\Z)$. It is not hard to see that the matrices $\pm E_\pi$ are conjugate to each other via reflection with respect to the edge containing the fixed point while $\pm E_{\frac{2\pi}{3}}$ are not, simply by the fact that they have different traces. Hence, there are three conjugacy classes (represented by $E_\pi, \pm E_{\frac{2\pi}{3}}$) in $GL(2,\Z)$ of elliptic matrices in $SL(2,\Z)$.

\textbf{Parabolic matrices:}  The fixed point of the action of a parabolic matrix in $PSL(2,\Z)$ is rational; thus, it is a common vertex of an infinite set of the triangles of $\D_F$. Since $PGL(2,\Z)$ acts transitively on the rational points, it is
not restrictive to assume that the fixed point of the translation is 0.

     \begin{figure}[ht]
   \begin{center}
      \includegraphics[ width=4.8cm, scale=0.3]{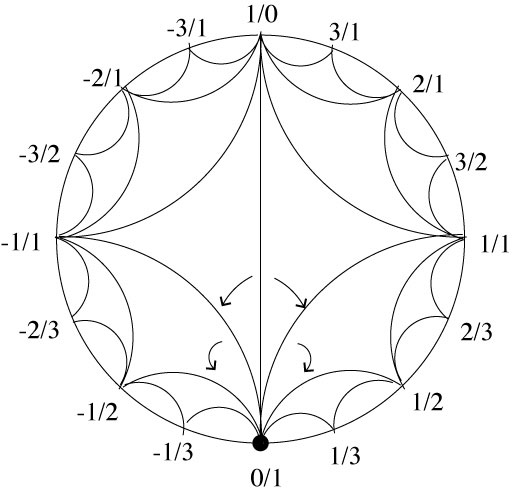}
       \caption{\small{Modular actions of parabolic matrices $P_n$.}}
       \label{parabolic}
      \end{center}
 \end{figure}

Hence, a parabolic element can shift the triangle $\Gamma$ by  an arbitrary number
$n$ of triangles to the right or to the left (see Fig.~\ref{parabolic})
fixing the point 0. The left shift is conjugated to the right shift by the
reflection with respect to the vertical line. Hence, the
equivalence classes in $PGL(2,\Z)$ are determined by the number
$n$ of shifts.
Such a shift can be represented by the matrix $P_n= \tiny{\left(\begin{array}{cc}
           1 & 0\\
           n & 1
         \end{array} \right)}$,  $n \in \N$.

The matrix $P_n\in PSL(2,\Z)$ corresponds to two matrices $\pm P_n$ in $SL(2,\Z)$. Having traces with opposite signs, they belong to two different conjugacy classes which are, thus, determined by the integer $\pm n$. Representatives of the conjugacy classes can be chosen as $\pm \tiny{\left(\begin{array}{cc}
           1 & 0\\
           n & 1
         \end{array} \right)}$, $n \in \N$.

\section{Conjugacy classes of hyperbolic matrices}\label{hiperbol}
\begin{figure}[ht]
   \begin{center}
      \includegraphics[ width=4.8cm, scale=0.3]{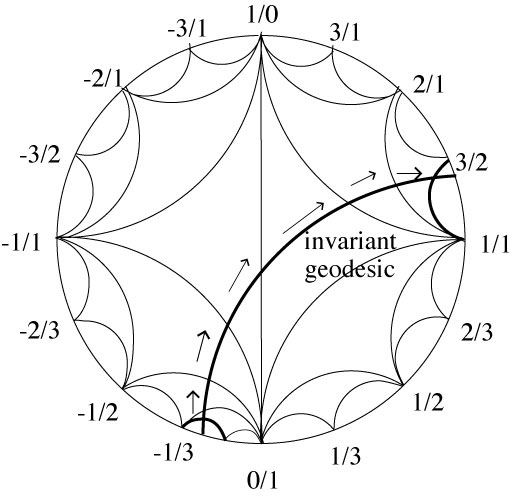}
       \caption{\small{Modular action of a hyperbolic matrix.}}
       \label{hyperbolic}
  \end{center}
\end{figure}

A hyperbolic matrix $A \in PSL(2,\Z)$ acts on $\D_F$ as translation fixing two irrational points. The geodesic (a semicircle), $l_A$, connecting these fixed points, oriented in the direction of translation, remains invariant under the translation, so $A$ preserves also the set of the triangles of $\D_F$ which are cut by $l_A$, see Fig.\ref{hyperbolic}. 

\begin{figure}[ht]
   \begin{center}
      \includegraphics[width=4.8cm, scale=0.27]{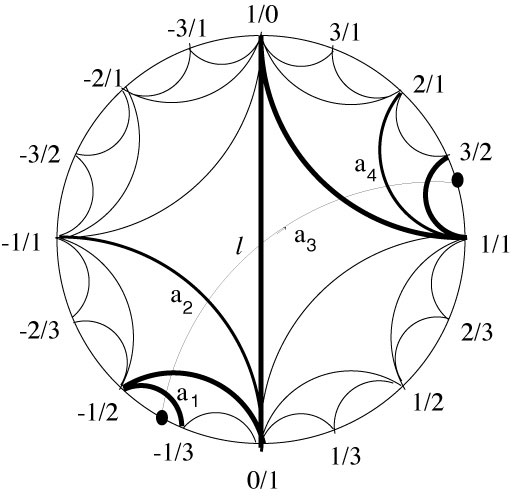}\hspace{.4cm}
       \caption{\small{Periodic pattern of the truncated triangles of the Farey tessellation.}}
       \label{period}
      \end{center}
\end{figure}

With respect to the orientation of $l_A$ (orientation is defined by the action of $A$),  such triangles are situated in two different ways: a set of triangles with a common vertex lying on the left of $l_A$ followed by a set of triangles with common vertex lying on the
right of $l_A$, see Fig.~\ref{period}.

Let us label right and left triangles by $R$ and  $L$,
respectively. Then we encode the arrangement of left and right
triangles with respect to $l_A$ as an infinite word, $\ldots LL
\ldots L RR \ldots RLL \ldots L \ldots$, of  2 letters. This word
is called  the {\em cutting word} of $l_A$. Let us fix a point $p$
at the intersection of $l_A$ with an edge of a triangle. Relative to this point, we obtain a sequence, $(a_1,
a_2, a_3,\ldots)_p$,  from the cutting word  where $a_{2i-1}$
stands for the number of consecutive triangles of one type while
$a_{2i}$ is the number of consecutive triangles of
the other type. For example, if the cutting word with respect to
$p$ reduced to the word $\underbrace{LL\ldots
L}_{a_1}\underbrace{RR\ldots R}_{a_2}\underbrace{LL\ldots
L}_{a_3}\ldots =L^{a_1}R^{a_2}L^{a_3}\ldots$, then we obtain
$(a_1, a_2,\ldots)_p$. This sequence is called the {\em cutting
sequence} relative to the point $p$.

Left and right triangles form a periodic pattern, and the action of $A$ is a shift by a period. Since the choice of $p$ is arbitrary the cutting sequence relative to $p$ is periodic after possibly some finite terms. Moreover, its period is of even length. Note that the choice of the point $p$ is not canonical; hence, we can encode the period only as a cycle, $[a_1 a_2\ldots a_{2n-1} a_{2n}]_A$, which we call the {\em cutting period-cycle associated to the matrix $A$}.

Because of the fact that $PGL(2, \Z)$ is the full symmetry group
of $\D_F$, the cutting period-cycle of a 
 hyperbolic matrix $A\in PSL(2,\Z)$ gives the complete invariant of the conjugacy class in
$PGL(2,\Z)$ of $A$. In other words, two matrices $A, B \in
PSL(2,\Z)$ are in the same conjugacy class in $PGL(2,\Z)$ if and
only if $[a_1 a_2\ldots a_{2n}]_A=[a_{\sigma(1)}
a_{\sigma(2)}\ldots a_{\sigma(2n)}]_B$ for a cyclic permutation
$\sigma$. Hence, we will denote the conjugacy classes in
$PGL(2,\Z)$ of hyperbolic matrices of $PSL(2,\Z)$ by the cycle
$[a_1a_2\ldots a_{2n}]$ (defined up to cyclic ordering).

It can be seen geometrically that with respect to the triangle $\Gamma$ a matrix representing a translation corresponding to the cutting period-cycle $[a_1 a_2 \ldots a_n]$ can be chosen as the following product of parabolic matrices. 

$$\small{\left(\begin{array}{cc}
           1 & a_1\\
           0 & 1
         \end{array} \right)
         \left(\begin{array}{cc}
           1 & 0\\
         a_2 & 1
         \end{array} \right)\cdots\left(\begin{array}{cc}
           1 & a_{2n-1}\\
           0 & 1
         \end{array} \right) \left(\begin{array}{cc}
           1 & 0\\
         a_{2n} & 1
         \end{array} \right)}.$$

For the sake of simplicity, let us
denote $U=\tiny{\left(\begin{array}{cc}
           1 & 1\\
           0 & 1
         \end{array} \right)}$ and $V=\tiny{
         \left(\begin{array}{cc}
           1 & 0\\
          1 & 1
         \end{array} \right)}$. Then the above product is written as $U^{a_1}V^{a_2}\ldots V^{a_{2n}}$. Note that $U$ is conjugate to $V$ in $PGL(2,\Z)$ but not in $PSL(2,\Z)$.

Let us note that in certain cases, namely if $l_A$ intersects the vertical line of $\D_F$ (since the action of $PGL(2,\Z)$ is transitive on the geodesics of $\D_F$, up to conjugation this property is always satisfied),  the cutting sequence of $l_A$ with respect to the point of intersection of $l_A$ with the vertical line is related to the continued fraction expansion of the fixed point $\xi$ which is the ``end point'' of $l_A$ with respect to the orientation. The corresponding theorem is due to C.~Series \cite{S2, S1}.

\begin{t.}\cite{S2, S1} Let $x > 1$, and let $l$ be any geodesic ray joining some point $p$
on the vertical line of $\D_F$ to $x$, oriented from $p$ to
$x$. Suppose that the cutting word of $l$ with respect to $p$ is
$L^{a_1}R^{a_2}L^{a_3}\ldots$. Then
$x=a_1+\frac{1}{a_2+\frac{1}{a_3+\cdots}}$. 

Note that if $0< x< 1$, then the sequence starts with $R$ and
$x=\frac{1}{a_1+\frac{1}{a_2+\frac{1}{a_3+\cdots}}}$. \\

If $x<0$ everything applies with $x$ replaced by $-x$ and with $R$
and $L$ interchanged.
\end{t.}

A matrix $A\in PSL(2,\Z)$ corresponds to $\pm A$ in $SL(2,\Z)$. Since $\pm A$ have
traces with opposite signs, the cutting period-cycle $[a_1 a_2 \ldots
a_{2n}]$, together with the sign of the trace determine the conjugacy classes of
$\pm A$ in $GL(2,\Z)$.  Representatives of them
can be chosen as $\pm U^{a_1} V^{a_2}\ldots
V^{a_{2n}}$.

\section{Real factorization of elliptic and parabolic matrices}\label{monrelf}

\begin{figure}[ht]
   \begin{center}
         \includegraphics[width=4.8cm, scale=0.3]{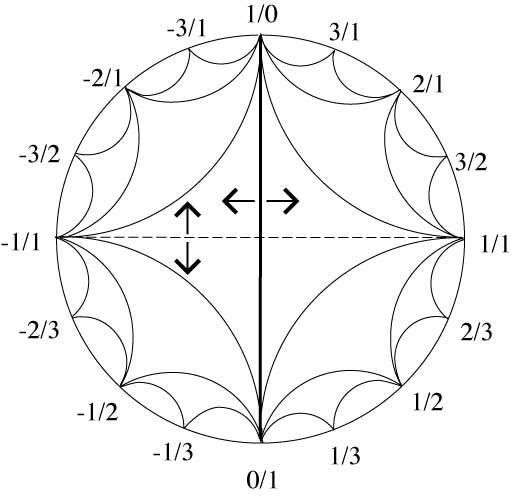}
       \caption{\small{Modular actions of the linear real structures.}}
       \label{yansimalar}
      \end{center}
 \end{figure}

The modular action of the linear real structures 
$\tiny{\left(\begin{array}{cc}
          1 & 0\\
          0 & -1
         \end{array} \right)},\tiny{\left(\begin{array}{cc}
           0 & 1\\
           1 & 0
         \end{array} \right)}$
on the hyperbolic plane $\D_{F}$ is $z\mapsto -\bar z$ and
$z\mapsto\frac1{\bar z}$, respectively. 
Geometrically, these are the reflections with respect to the vertical and respectively the
horizontal lines, see Fig.~\ref{yansimalar}. In particular, the first
reflection takes the basic triangle $\Gamma$  to the triangle with vertices $\{0, -1, \infty \}$,
and the second one takes $\Gamma$ to itself.

\begin{t.}\label{elpar}
All elliptic and parabolic matrices in $SL(2,\Z)$ are products of two linear real structures.
\end{t.}

\noindent {\it Proof:} 
The explicit real decompositions for the representatives of  the conjugacy classes of
the elliptic matrices are given below:

$$\begin{array}{rcl}
E_{\frac{2\pi}{3}}&=&\small{\left(\begin{array}{cc}
           0 & 1\\
           -1 & 1
         \end{array} \right)=\left(\begin{array}{cc}
           1 & 0\\
           1 & -1
         \end{array} \right)\left(\begin{array}{cc}
           0 & 1\\
           1 & 0
         \end{array} \right)}\\
-E_{\frac{2\pi}{3}}&\cong &\small{\left(\begin{array}{cc}
           -1 & 1\\
           -1 & 0
         \end{array} \right)=\left(\begin{array}{cc}
           1 & -1\\
           0 & -1
         \end{array} \right)\left(\begin{array}{cc}
           0 & 1\\
           1 & 0
         \end{array} \right)}\\
 E_\pi&=&\small{\left(\begin{array}{cc}
           0 & 1\\
          -1 & 0
         \end{array} \right)=\left(\begin{array}{cc}
           1 & 0\\
           0 & -1
         \end{array} \right)\left(\begin{array}{cc}
           0 & 1\\
           1 & 0
        \end{array} \right)}.\end{array}$$

\begin{figure}[ht]
   \begin{center}
         \includegraphics[width=4cm, scale=0.3]{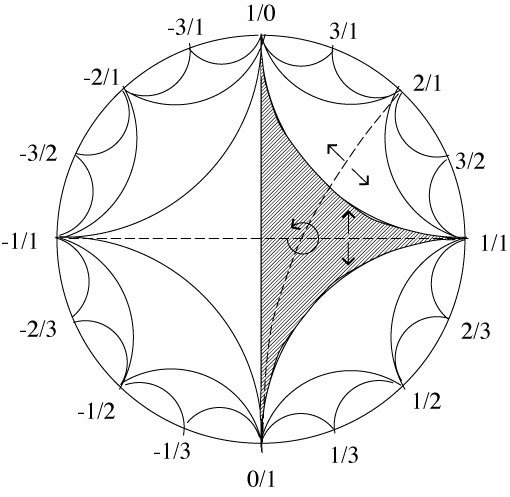}
      \includegraphics[width=4cm, scale=0.3]{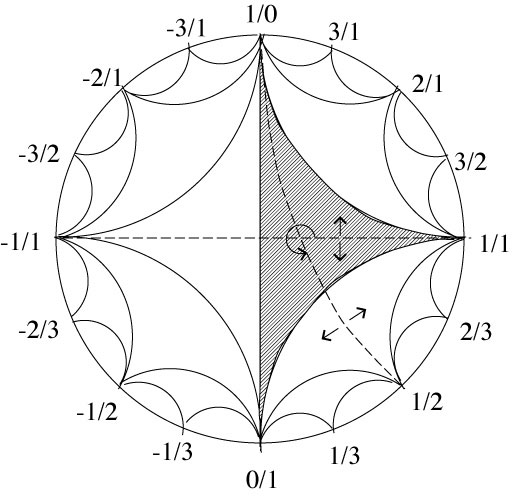}
      \includegraphics[width=4cm, scale=0.3]{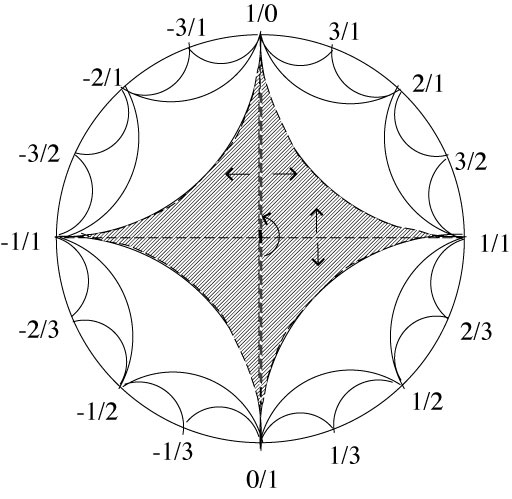}
       \caption{\small{Decompositions of the modular actions of the elliptic matrices.}}
       \label{eref}
      \end{center}
 \end{figure}

Fig.~\ref{eref} illustrates geometrically the above decompositions in terms of the
corresponding modular action of the matrices.

Real decompositions for  the representatives of the conjugacy classes of the parabolic matrices are as
follows:

         $$P_n=\small{\left(\begin{array}{cc}
           1 & 0\\
           n & 1
         \end{array} \right)=\left(\begin{array}{cc}
           1 & 0\\
           n & -1
         \end{array} \right)\left(\begin{array}{cc}
           1 & 0\\
           0 & -1
         \end{array} \right)}$$

           $$-P_n=\small{\left(\begin{array}{cc}
           -1 & 0\\
           -n & -1
         \end{array} \right)=\left(\begin{array}{cc}
           1 & 0\\
           n & -1
         \end{array} \right)\left(\begin{array}{cc}
           -1 & 0\\
            0 & 1
         \end{array} \right)}.$$
\hfill  $\Box$\\

\begin{figure}[ht]
   \begin{center}
      \includegraphics[width=4cm, scale=0.3]{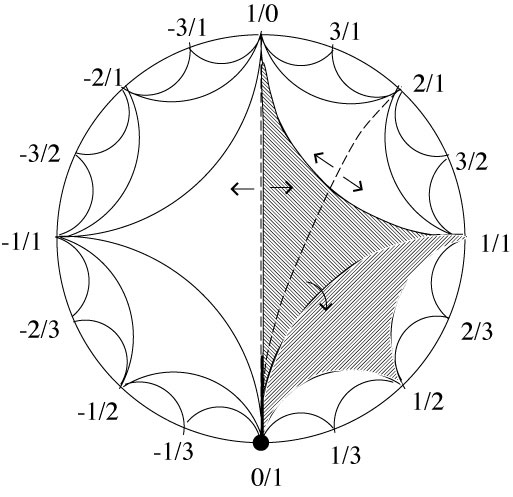}\hspace{.5cm}
      \includegraphics[width=4cm, scale=0.3]{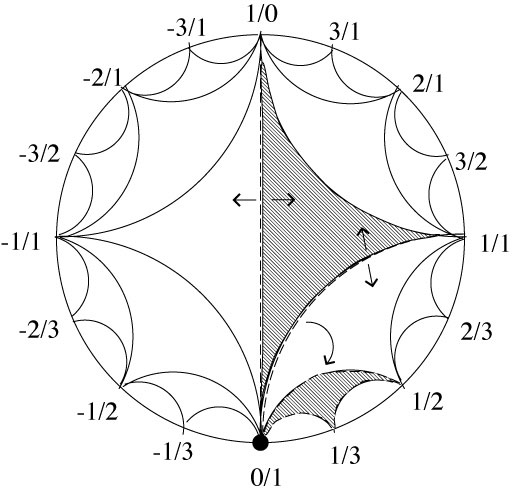}
       \caption{\small{Decompositions of the modular actions of the parabolic matrices $P_1,\, P_2$.}}
       \label{p}
      \end{center}
 \end{figure}

Fig.~\ref{p} shows  the real decompositions of the
modular action of  matrices $\tiny{\left(\begin{array}{cc}
           1 & 0\\
           n & 1
         \end{array} \right)}$ $n=1,2$.

\section{The criterion of factorizability for hyperbolic matrices}

\begin{l.} \label{2sym} If the cutting period-cycle of a hyperbolic matrix
$A$ is $[a_1 a_2 \ldots a_{2n}]$, then the cutting period-cycle of
$A^{-1}$ is $[a_{2n} a_{2n-1} \ldots a_{1}]$.
\end{l.}

\noindent {\it Proof:} 
Note that $l_A= l_{A^{-1}}$ with the opposite orientation. So the cutting word of $A^{-1}$ can be obtained from the cutting word of $A$ by taking the mirror image of the word and interchanging $L$ with $R$. Interchanging $L$ and $R$ does not effect the cutting period-cycle; hence, the cutting period-cycle of $A^{-1}$ is the reverse $[a_{2n} a_{2n-1}\ldots a_1]$ of the cutting period-cycle $[a_{1} a_{2}\ldots a_{2n}]$ of $A$.
\hfill  $\Box$\\

\begin{d.}
A finite sequence $(a_1a_2\ldots a_k)$ is called \emph{palindromic} if it is equal to the reversed sequence $(a_k a_{k-1}\ldots a_1)$.
We call $k$ the \emph{length of the sequence}. 

A finite sequence is called \emph{bipalindromic} if it can be subdivided into two palindromic sequences. 
\end{d.}

\begin{d.}
A cutting period-cycle is called \emph{bipalindromic} if there is a cyclic permutation such that the permuted period is bipalindromic.

In particular, if the cutting period-cycle is subdivided into two palindromic sequences of odd length (respectively even length) we call it \emph{odd-bipalindromic} (respectively \emph{even-palindromic}). 
\end{d.}

To illustrate, the period $[1213]$ is odd-bipalindromic, while the period $[1122]$ is even-bipalindromic. 

\begin{l.} If $A^{-1}=Q^{-1}AQ$ for some $Q\in PGL(2,\Z)$, then the cutting period-cycle $[a_1 a_2\ldots a_{2n}]_{A}$ is bipalindromic. 
\end{l.}

\noindent {\it Proof:}  If  $A$ and $A^{-1}$ are in the same conjugacy class in $PGL(2,\Z)$, then they have the same cutting period-cycles up to cyclic permutation. By Lemma~\ref{2sym}  we know that the cutting period-cycle of $A^{-1}$ is $[a_{2n} a_{2n-1}\ldots a_{1}]$ while the cutting period-cycle of $A$ is $[a_1 a_2\ldots a_{2n}]$. Hence, $[a_{\sigma(1)}  a_{\sigma(2)} \ldots a_{\sigma{(2n)}}]= [a_{2n} a_{2n-1} \ldots a_{1}]$ for some cyclic permutation $\sigma$.  
Without loss of generality, let us assume that $\sigma$ is a shift by $k=2n-l$, so we have $(a_{l+1}  a_{l+2} \ldots a_{2n} a_{1} a_{2} \ldots a_{l})= (a_{2n} a_{2n-1} \ldots a_{l+1}a_{l}a_{l-1} \ldots a_{1})$ as finite sequences. Thus,  $(a_{l+1}  a_{l+1} \ldots a_{2n})= (a_{2n} a_{2n-1} \ldots a_{l+1})$ and  $(a_{1} a_{2} \ldots a_{l})=(a_{l} a_{l-1} \ldots a_{1})$ which implies that the cutting period-cycle is bipalindromic.
\hfill  $\Box$\\

Note that if the cutting period-cycle is odd-bipalindromic, then the symmetry of the palindromic pieces lifts to a symmetry of the left and the right triangles corresponding to the cutting period-cycle. However, this is not true for even-bipalindromic periods. In the case of $[1213]$, we have $121\sim LR^2L=LRRL$ and $3\sim R^3=RRR$, whereas of $[1122]$, we have $11\sim LR$ and $22\sim L^2R^2=LLRR$.

\begin{t.} \label{hiper}
A hyperbolic matrix $A$ is a product of two linear real structures if and only if its cutting period-cycle  $[a_1 a_2 \ldots a_{2n}]_A$ is odd-bipalindromic.
\end{t.}

\begin{l.}\label{c(l)=l} Let $A \in PSL(2,\Z)$ such that $A^{-1}=Q^{-1}AQ$ for some $Q \in PGL(2,\Z)$, and let $l_A$ be the geodesic invariant under the action of $A$. Then $Q(l_A)=l_A$.
\end{l.}

\noindent {\it Proof:} 
Clearly, if $A(l_A)=l_A$, then $A^{-1}(l_A)=l_A$. Hence, $$A^{-1}(l_A)=Q^{-1}AQ(l_A) \Leftrightarrow  Q(l_A)=A(Q(l_A)).$$
By the uniqueness of the invariant geodesic we obtain $Q(l_A)=l_A$.
\hfill  $\Box$\\

\begin{l.}\label{orpres}
Let $A, Q, l_A$ be as above. If the cutting period-cycle  $[a_1
a_2\ldots a_{2n}]_A$ is even-bipalindromic, then $Q$ is
orientation preserving.
\end{l.}

\noindent {\it Proof:} 
By means of Lemma~\ref{c(l)=l}, we have $Q(l_A)=l_A$; hence, $Q$  preserves triangles meeting $l_A$.  The action of $Q$ on $\D_F$ is a linear fractional transformation; thus, it preserves the angles. An analysis on the angles at the meeting points of $l_A$ and the edges of the triangles will forbid the existence of the orientation reversing map when the cutting period-cycle is even-bipalindromic.

\begin{figure}[ht]
   \begin{center}
   \includegraphics[width=4.8cm, scale=0.3]{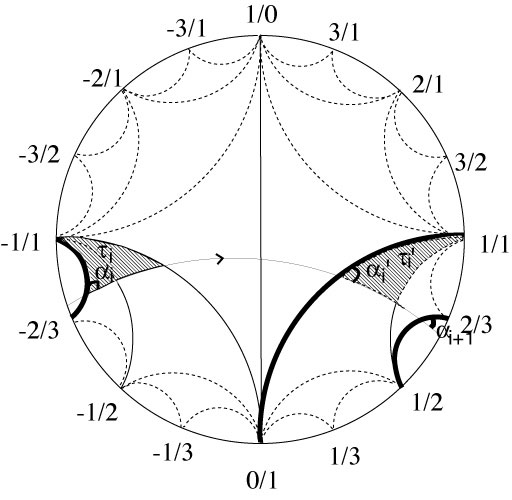}\hspace{.5cm}
\includegraphics[width=6.5cm, scale=0.3] {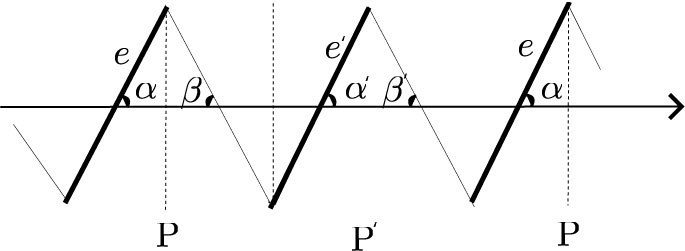}
       \caption{Triangles $\tau_{i}$ and $\tau'_{i}$ and the interior angles $\alpha_{i}, \alpha'_{i}$, $\beta_{i}$, $\beta'_{i}$.}
       \label{acilar}
  \end{center}
\end{figure}

Let us assume that the cutting period-cycle has the form:
 $$[\underbrace{a_1 a_2\ldots a_k a_k \ldots a_2 a_1}_{P}\underbrace{a'_1 a'_2 \ldots a'_s a'_s \ldots a_2 a'_1}_{P'}]$$ where $s+k=n$ and $P$ and $P'$ are two palindromic pieces.
Substituting the pieces $P, P'$ to the cutting sequence, we obtain a sequence of  $P$ and $P'$ of the form $\ldots PP'PP'\ldots$.  Clearly, the action of the matrix $A$ corresponds to a shift by two: it  takes $P$ to $P$, $P'$ to $P'$. 
Let us consider the edges which separate the triangles corresponding to $P$ from the triangles corresponding to $P'$. They are of two types:  with respect to the orientation of  $l_{A}$, we encounter the edges where we move from $P$ to $P'$, and inversely, the edges where we pass from $P'$ to $P$. We denote such edges by $e_{i}$ and $e'_{i}$ respectively, see Fig.~\ref{acilar}.  (We enumerate the edges with respect to an auxiliary point $p$ fixed on $l_{A}$.)

 Each triangle of $\D_\F$ which is cut by $l_A$ splits into two pieces, one of which is a triangle. Let $\tau_i$ (respectively $\tau'_i$) denote the triangle having one edge $e_i$ (respectively $e'_i$) and obtained as the union of  triangle-pieces of the triangles of $\D_F$ with a common vertex on one side of $l_A$ (respectively on the other side of $l_{A}$), see Fig.~\ref{acilar}. Let $\alpha_i$ (respectively $\alpha'_i$) be the interior angles of $\tau_i$ (respectively $\tau'_i$) between the edges $e_i$ (respectively $e'_i$) and $l_A$.
In addition, let $\beta_i$ (respectively $\beta'_i$) be the other interior angle of $\tau_i$ (respectively $\tau'_i$) on $l_A$.

Note that since $A$ shifts triangles  by the period $PP'$, $A$ takes $\alpha_i$ to $\alpha_{i+1}$ (respectively $\alpha'_i$ to $\alpha'_{i+1}$). Hence, all $\alpha_i$ (respectively all $\alpha'_i$) are equal. Let $\alpha=\alpha_i$ for all $i$ (respectively  $\alpha'=\alpha'_i$ for all $i$).

\begin{figure}[ht]
   \begin{center}
   \includegraphics[width=7.5cm, scale=0.3]{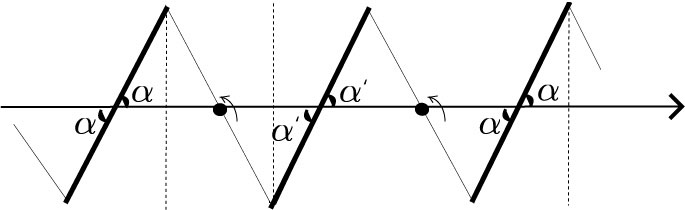}
       \caption{Elliptic rotations of even-bipalindromic cutting period-cycles.}
       \label{don}
  \end{center}
\end{figure}

The crucial observation is that if the cutting period cycle is even-bipalindromic, then there is an elliptic matrix in the conjugacy classes of $E_{\pi}$ which fixes the point of intersection of $l_A$ with the middle edge of $P$ or $P'$, see Fig.~\ref{don}. Such matrix interchanges the edges $e_i$ to $e'_i$. Hence, $\alpha=\alpha'$. 
(In the same way, we obtain $\beta=\beta_i=\beta'_i$ for all $i$.)

Let us assume that $\alpha< \frac{\pi}{2}$ (otherwise, we can replace  $\alpha$ with $\beta$.)
We choose an orientation of $\D_F$ by specifying $(v_1, v_2)$ where $v_1$ is a tangent vector of $l_A$ and $v_2$ is the tangent vector of $e_i$ (or $e'_i$) such that the angle $\alpha$ between $v_1$ and $v_2$ is $\alpha < \frac{\pi}{2}$.  The matrix  $Q$ takes $(v_1, v_2)$ to itself since it preserves $l_A$ and the set of edges $e_{i}, e'_{i}$; hence, it preserves the angles between the two. However, an orientation reversing map can not preserve both the angle  $\alpha < \frac{\pi}{2}$ between the vectors $(v_1, v_2)$ and the vectors $v_{1},v_{2}$ at the same time. Thus, $Q$ is orientation preserving.\hfill  $\Box$\\

\noindent {\it Proof of Theorem \ref{hiper}:}
If $A$ is  a product of two linear real
structures. Then the cutting period-cycle is odd-bipalindromic by
Lemma~\ref{orpres}.

If the cutting period-cycle is odd-bipalindromic, then up
to cyclic ordering, it has two palindromic
pieces of odd length. Let us assume that the cutting period-cycle
is of the form $$\pm[a_1 a_2\ldots a_k a_{k+1} a_k\ldots a_2 a_1 a'_1
a'_2\ldots a'_s a'_{s+1} a'_s\ldots a'_2 a'_1]$$ where
$(2k+1)+(2s+1)=2n$. Then for some $R\in PGL(2,\Z)$, we have
$B=R^{-1}AR$ such that $B=\pm U^{a_1}V^{a_2}\ldots U^{a_2}
V^{a_1}U^{a'_1}V^{a'_2}\ldots U^{a'_2} V^{a'_1}.$ Matrices
$U^{a_i}$ and $V^{a_i}$ have the following real decompositions:
$$ {U^{a_i}=\small{\left(\begin{array}{cc}
           1 & -a_i\\
           0 & -1
         \end{array} \right)\left(\begin{array}{cc}
           1 & 0\\
           0& -1
         \end{array} \right)}} \, \textrm{and} \, {V^{a_i}=\small{\left(\begin{array}{cc}
           1 & 0\\
           0 & -1
         \end{array} \right)\left(\begin{array}{cc}
           1 & 0\\
          -a_i& -1
         \end{array} \right)}}.$$
         
\vspace{0.3cm}

Hence, the product $U^{a_1}V^{a_2}\ldots U^{a_2}V^{a_1}$ can be
rewritten as follows:

 $$\small{\left(\begin{array}{cc}
           1 & -a_1\\
           0& -1
         \end{array} \right)\cdots\left(\begin{array}{cc}
           1 & 0\\
           -a_k& -1
         \end{array} \right)\left(\begin{array}{cc}
           1 & -a_{k+1}\\
           0& -1
         \end{array} \right)\left(\begin{array}{cc}
           1 & 0\\
           -a_k& -1
         \end{array} \right)\cdots\left(\begin{array}{cc}
           1 & -a_1\\
           0 & -1
         \end{array} \right)}. $$
         
      \vspace{0.3cm}   
This gives a linear real structure, since it is a conjugate of
$\tiny{\left(\begin{array}{cc}
           1 & -a_{k+1}\\
           0& -1
         \end{array} \right)}$. Similarly, the product $U^{a'_1}V^{a'_2}\ldots U^{a'_2}V^{a'_1}$ gives a linear real structure conjugate to $\tiny{\left(\begin{array}{cc}
           1 & -a_{s+1}\\
           0& -1
         \end{array} \right)}$. 
         \hfill  $\Box$\\

\begin{t.} \label{wrr}
A matrix $A \in SL(2,\Z)$ is a product of two linear real structures if and only if  there is a matrix $Q$ with $\det Q=-1$ such that $A^{-1}=Q^{-1}AQ$.
\end{t.}

\noindent {\it Proof:}  Necessity of the condition is trivial.
As for the converse, we only need to consider the case of hyperbolic matrices. 
Let $A$ be a hyperbolic matrix such that $ A^{-1}=Q^{-1} A Q$, for some $Q$ with $\det{Q}=-1$.
Then the cutting period-cycle $[a_1a_2\ldots a_{2n-1} a_{2n}]_{A}$ is odd-bipalindromic by Lemma~\ref{orpres}. So  $A$ is real by
Theorem~\ref{hiper}.
\hfill  $\Box$\\

Theorem~\ref{wrr} together with Proposition~\ref{s1fcc}  lead to the following corollary. 

\begin{c.}
An elliptic $F$-fibration is real if and only if it is weakly real.
 \hfill  $\Box$\\
\end{c.}

\end{document}